\documentclass[10pt,twocolumn]{article}
\usepackage[english]{babel}
\usepackage{amsfonts}
\usepackage[dvips]{graphicx}

\begin{document}

\title{\bf Integrability of Newton ovals, computation of air damper inlets}

\date{2011, September}

\author{\bf Gianluca Argentini \\
\normalsize{Research \& Development Dept., Riello Burners - Italy}\\
\normalsize gianluca.argentini@rielloburners.com \\
\normalsize gianluca.argentini@gmail.com \\}

\maketitle

\noindent{\bf Abstract}\\
About global and local algebraic integrability of ovals. A contribution to clarify Newton results and relative comments on his work done by Arnol'd and Pourciau. A possibile application to air damper sections computation is offered, as example of unexpected link between pure mathematics and industrial technology.\\


Following Arnol'd (\cite{arnold}), a curve in the real plane is said to be {\it algebraic} if its points $(x,y)$ satisfy an equation $P(x,y)=0$ where $P$ is a non-zero polynomial.
Following Newton in its {\it Principia} (\cite{chandra}), an {\it oval} is a closed convex algebraic curve; following Newton (and Chandrasekhar) again, an oval is said to be {\it smooth} if at every point there is the tangent line and its position varies with continuity when the point varies with continuity along the curve.

\begin{figure}[ht!]
	\begin{center}
	\includegraphics[width=7cm]{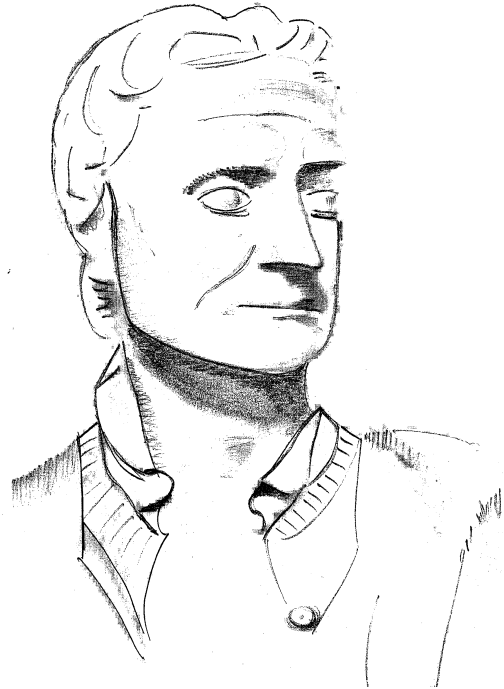}
	\caption{{\it Newton} {\small (pencil drawing by G.Argentini, 2011)}.}
	\label{Newton}
	\end{center}
\end{figure}

{\noindent} An oval is said to be {\it algebraically squarable} if exists a polynomial $Q(p,q,r,s)$ such that the area $S$ of an arbitrary segment formed by a cartesian line $ax+by+c=0$ satisfies the equation $Q(S,a,b,c)=0$ (\cite{arnold},\cite{pourciau}). Newton, during its studies on planets trajectories, has proven in its Lemma 28 of the First Book of Principia the following theorem (see (\cite{arnold}) for the original Newton proof, which is based on a geometric and algebric, not analytic, argument):\\

{\it If an oval is smooth, it is not algebraically squarable}.\\

\begin{figure}[ht!]
	\begin{center}
	\includegraphics[width=7cm]{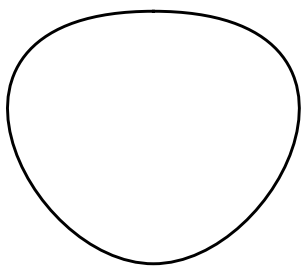}
	\caption{{\it The convex Newton apple}.}
	\label{appleNewton}
	\end{center}
\end{figure}

{\noindent} For example, Kepler elliptical ovals are not algebraically squarable, and so the convex Newton apple (Fig.\ref{appleNewton}; see forward for its equation). The unit square of equation $x^2y^2-x^2y-xy^2+xy=0$ is algebraically squarable and a polynomial satisfying the definition is $Q(S,m,q)=4m^2S^2-4m^2S-2mq+q^2+2mq^2-2q^3+q^4$.\\

\begin{figure}[ht!]
	\begin{center}
	\includegraphics[width=7cm]{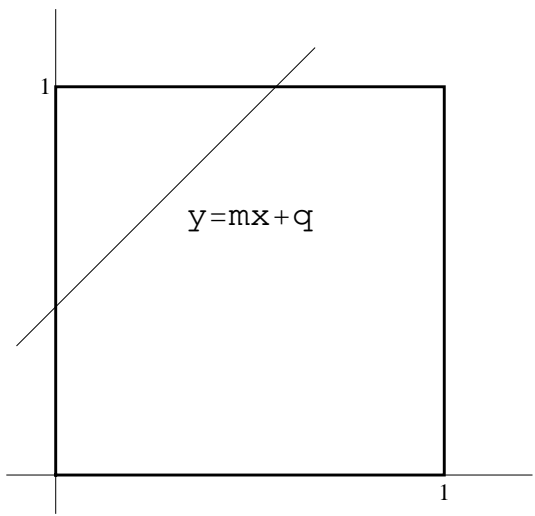}
	\caption{{\it The squarable square}.}
	\label{square}
	\end{center}
\end{figure}

\noindent In their works, Arnol'd and Pourciau describe how the Newton result about ovals integrability can be interpreted in a weaker sense too, using the same logical argumentation used by Newton in his proof. They introduce the concept of {\it locally algebraically integrability}, without giving a clear formal definition. But from the examples and from the proofs which they offer, one can deduce the following definition:\\

{\it An oval is said to be locally algebraically squarable if for each point $P_0$ exist a line $L_0$ of equation $a_0x+b_0y+c_0=0$ and a polynomial $Q_0(p,q,r,s)$ such that the areas $S$ of the segments formed by all the lines $ax+by+c=0$, in a suitable geometrical neighbourhood of the line $L_0$, satisfy the equation $Q_0(S,a,b,c)=0$}.\\

\noindent In practice, in (global) algebraic integrability for the whole oval there is a polynomial $Q$ such that $Q(S,a,b,c)=0$, while in local algebraic integrability for each point $P_0$ of the oval there is a polynomial $Q_0$ such that $Q_0(S,a,b,c)=0$.

As Newton understood (see \cite{arnold}), global and local algebraic integrability seem to depend on the degree of smoothness of the curve.
A curve $F(x,y)$ is said to be of class $C^m$ at its point $P_0=(x_0,y_0)$ if its partial derivatives until the $m$-th order exist at $P_0$, they are continuous, and the gradient $(\partial_xF(x_0,y_0), \partial_yF(x_0,y_0))$ is not the null vector (that is $P_0$ is {\it regular}). A curve is said to be of class $C^m$ if all its points are regular and its partial derivatives until the $m$-th order exist and are continuous. If also the curve is infinitely derivable with continuity, it is of class $C^\infty$. The curve is said to be {\it analytic} at $P_0$ if it is expandible in a convergent power series of integer exponents in a neighbourhood of $P_0.$\\

{\bf Theorem} {\it An algebraic curve of class $C^1$ at its point $P_0$, is analytic at $P_0$}.\\

{\it Proof}. Suppose $\partial_y F(x_0,y_0) \neq 0$. From the theorem of implicit function, there is a function $y=\phi(x)$ such that $F(x,\phi(x))=0$ for all the $x$ in a neighbourhood of $P_0$. Also, $\phi$ is derivable and its first derivative is

\begin{equation}\label{phi_first}
	\phi'(x) = - \partial_x F(x,\phi(x))/\partial_y F(x,\phi(x))
\end{equation}

Note that, using permanence of sign, the neighbourhood of $P_0$ can be chosen so that $\partial_y F(x,\phi(x)) > 0$ on it. Therefore, being $F(x,y)$ partially derivable indefinitely, from ($\ref{phi_first}$) the successive derivatives $\phi''$, $\phi'''$, ... have only integer powers of $\partial_y F(x,\phi(x))$ at denominator, and we deduce that $\phi$ is indefinitely derivable. So locally the curve $F(x,y)$ is the graph of an infinitely differentiable function $y=\phi(x)$. From a theorem of Newton (\cite{arnold}), it follows that $F$ is analytic at $P_0$. $\square$\\

A curve is {\it analytic} if it is analytic at every point. Newton has proven (\cite{arnold}) that a $C^\infty$ oval is analytic and it is not algebraically squarable, even locally. From previous theorem it can be deduced that\\

{\it A $C^1$ oval is not algebraically squarable, even locally}.\\

Hence, if you want an oval locally algebraically squarable, you must search among ovals with at least a {\it singular} point, that is a point $P_0=(x_0,y_0)$ such that $\partial_xF(x_0,y_0)=\partial_yF(x_0,y_0)=0$. Note that an oval can be smooth in the sense of Newton, that is to have continuous tangent line at every point, but it could not be of class $C^1$ (in the paper of Pourciau, smooth is equivalent to $C^1$, but this seems not to agree with Newton assumptions). Arnol'd (\cite{arnold}) offers the example of the parametric curve $x=(t^2-1)^2, \hspace{0.2cm} y=t^3-t$ for $-1 \leq t \leq 1$, which have continuous and never null tangent vector $(4t(t^2-1), 3t^2-1)$, that is the curve is smooth. But if we consider its algebraic equation, that is $y^4-2xy^2-x^3+x^2=0$, we can see that the origin $(0,0)$ is the unique singular point, that is the curve is not of class $C^1$. Therefore, from Newton this oval is not algebraic squarable. It is locally algebraic squarable (see (\cite{arnold}), but it is not a counterexample to Newton Lemma 28 as declared in (\cite{pourciau}), because Newton denied local integrability only for $C^1$ algebraic ovals.\\

Following Newton and previous theorem, a $C^1$ oval has at every point a local representation ($y=y(x)$ or $x=x(y)$) as an integer power series. If an oval has a singular point, it has a local representation as a fractional power series of the kind $y=a_0+a_1x^\frac{1}{N}+a_2x^\frac{2}{N}+...$ for an integer $N \geq 2$, as shown by Puiseux (\cite{wall}). Note that Newton knew how to expand a curve in a fractional power series using the method today known as {\it Newton polygon} (\cite{arnold},\cite{wall}).\\

\begin{figure}[ht!]
	\begin{center}
	\includegraphics[width=7cm]{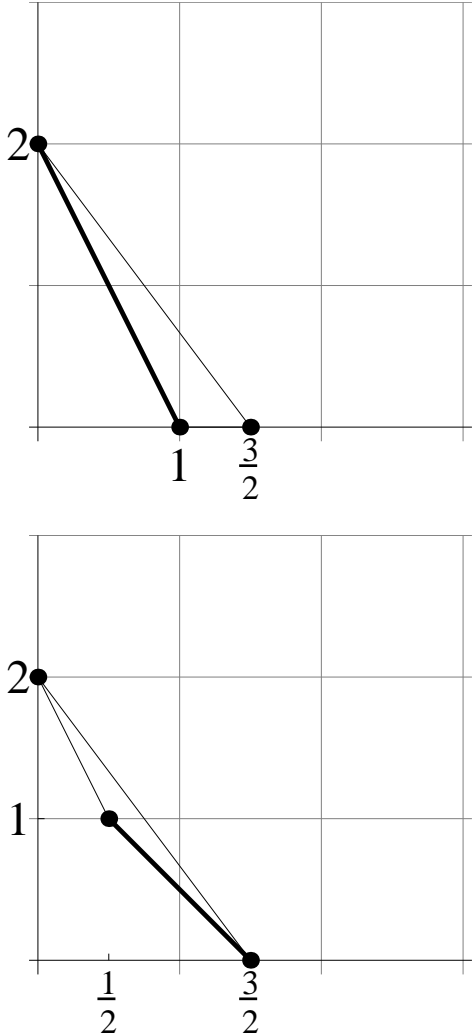}
	\caption{{\it Newton polygons for fractional power series expansion of the oval} $y^4-2xy^2-x^3+x^2=0$.}
	\label{NewtonPolygons}
	\end{center}
\end{figure}

Consider previous Arnol'd oval $y^4-2xy^2-x^3+x^2=(y^2-x)^2-x^3=0$. How we can find a local representation of kind $y=y(x)$ at the singular point $(0,0)$? At first step, from $y^2=x+x^\frac{3}{2}$ using Newton polygon (Fig.\ref{NewtonPolygons}) we have the equation $y^2-x=0$, from which follows (considering only the positive branch) $y=x^\frac{1}{2}$. Now let be $y=x^\frac{1}{2}+z$, and we compute $z$ using the new Newton polygon. Consider now the equation $2x^\frac{1}{2}z-x^\frac{3}{2}=0$, which gives $z=\frac{1}{2}x$. So we have the following local representation near the origin: $y=x^\frac{1}{2}+\frac{1}{2}x$. Note that this representation are the first two terms of the Puiseux series $y=\sum_{n=1}^{+\infty} a_n x^\frac{n}{2}$ (a {\it Mathematica} expansion gives $x^\frac{1}{2}+\frac{1}{2}x-\frac{1}{8}x^\frac{3}{2}+\frac{1}{16}x^2-\frac{5}{128}x^\frac{5}{2}+O[x]^3$). So we have a singular point where the oval is locally represented by a power series where not all the exponents are integer; the oval is not analytic but it is smooth and locally algebraically squarable.

\begin{figure}[ht!]
	\begin{center}
	\includegraphics[width=7cm]{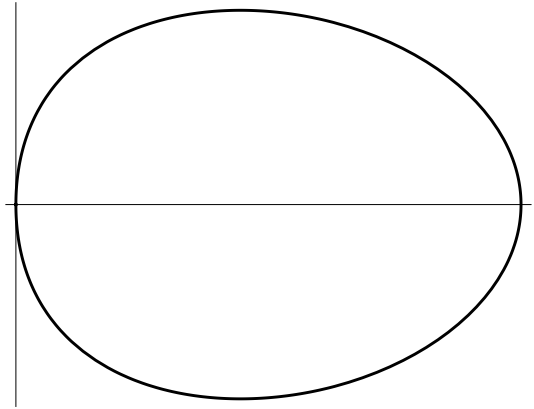}
	\caption{{\it The Arnol'd oval} $y^4-2xy^2-x^3+x^2=0$.}
	\label{algebraicCurve}
	\end{center}
\end{figure}

Let $(x_0,y_0)$ be a point of an oval. The oval is said to have a {\it $(x_0,y_0)$-centered parametrization} if there are two rational functions $g(t)$, $f(t)$ in the interval $[a,b]$, $a < b$, such that $x=g(t)$, $y=f(t)$ for every point of the oval, $g(a)=g(b)=x_0$, $f(a)=f(b)=y_0$, and $g(t) \neq x_0$ for every $t$ in $(a,b)$. Note that $(g(a),f(a))$ can be an angular or a smooth point.\\

{\bf Theorem} {\it An oval having a singular point $(x_0,y_0)$ and a $(x_0,y_0)$-centered parametrization, is locally algebraically squarable}.\\

{\it Proof}. We can suppose, without loss of generality, that the singular point $(x_0,y_0)$ is the origin $(0,0)$, because a translation $X=x-x_0$, $Y=y-y_0$ of cartesian axes doesn't modify the character of a singular point of a curve $F(x,y)=0$ (if $G(X,Y)=F(x(X),y(Y))$, then $\partial_XG=\partial_xF\partial_Xx=\partial_xF$ and $\partial_YG=\partial_yF$). Let $P_0$ a point on the oval distinct from the origin and $t_0$ the value of the parameter such that $(g(t_0),f(t_0))=P_0$. Note that, by definition of centered parametrization, $g(t_0) \neq 0$. Then the straight line from origin to $P_0$ has equation $y=mx$ with $m=f(t_0)/g(t_0)$. Note that $m$ has a rational dependence on $t_0$. Suppose that the oval is described clockwising by $t$; the area $S$ of the (upper) segment can be computed by integration of the $y$ values on $x$ minus the area $1/2g(t_0)f(t_0)$ of the triangle down the straight line. So

\begin{equation}
	S = \int_a^{t_0} f(t)g'(t)dt - \frac{1}{2}g(t_0)f(t_0)
\end{equation}

\noindent This expression of $S$ has a rational dependence on $t_0$. Eliminating $t_0$ from the expressions of $m$ and $S$ by computation of their resultant (\cite{cox}), one find an algebraic equation between $m$ and $S$, wich is valid for all the straight lines of the kind $y=nx$ with $n$ variable in a neighbourhood of $m$.\\
Let us consider now the origin. We can choose a segment delimited by a vertical line of kind $x=c$ whose intersections with the oval let be represented by the two parameter values $t_1$ for the upper and $t_2$ for the lower; again, $c$ has a rational dependence $f_1(t_1)$ on $t_1$ and a rational dependence $f_2(t_2)$ on $t_2$. The area of the segment is 

\begin{equation}
	S = \int_a^{t_1} f(t)g'(t)dt + \int_{t_2}^b f(t)g'(t)dt
\end{equation}

\noindent so again, eliminating $t_1$ and $t_2$ from the set of equations $S=S(t_1,t_2)$, $c=f_1(t_1)$, $c=f_2(t_2)$, we obtain an algebraic relation between $S$ and $c$. $\square$\\

\begin{figure}[ht!]
	\begin{center}
	\includegraphics[height=7cm]{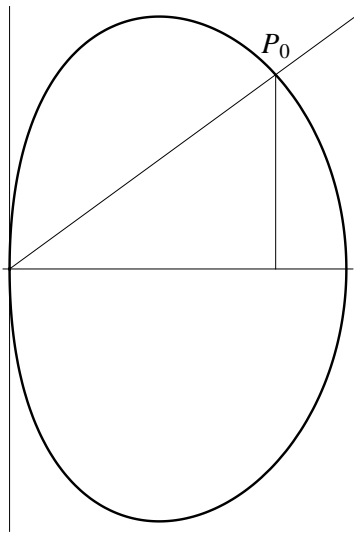}
	\caption{{\it Construction of the segments for local integrability}.}
	\label{segmentForProof}
	\end{center}
\end{figure}

\noindent {\it Remark}. From definition, the denominators of the rational functions $g$ and $f$ can't have zeros in the interval $(a,b)$. So at points distinct from $(x_0,y_0)$ the oval is smooth in the sense of Newton.\\

The hypothesis of singularity for the point $(x_0,y_0)$ seems to remain hidden in the proof. But if the oval has only regular points, it shall be a $C^1$ oval, so it coudn't be locally algebraically squarable.\\

\noindent {\it Example}. The couple of functions $x=(t^2-1)^2, \hspace{0.2cm} y=t^3-t$ is a $(0,0)$-centered parametrization for the Arnol'd oval.\\

\noindent A {\it Bezier oval} is a plane convex close smooth (in Newton sense) Bezier curve (\cite{salomon}) with $(0,0)$ as first and last control point. Note that for such a curve $(0,0)$ is a singular point and the usual Bernstein parametrization is $(0,0)$-centered. Therefore a Bezier oval is locally algebraically squarable.\\

\noindent {\it Example}. The convex Newton apple of Fig.\ref{appleNewton} is actually a Bezier oval wich can be constructed using the list of control points $\{(0, 0), (-a, 0), (-b, c), (0, c), (b, c), (a, 0), (0, 0)\}$, with $a$, $b$ and $c$ geometric parameters; its algebraic equation $F(x,y)=0$ has degree six. The origin is singular. It is locally algebraically squarable.\\

An air damper is a mechanical component which has the function of varying the air flow rate in fluid dynamical systems, such as ventilating structures. Usually there is an hole on a material boundary and a cover (mechanical doorlet or shutter or simply damper) manually or automatically moved for partially opening or closing the hole section.

\begin{figure}[ht!]
	\begin{center}
	\includegraphics[height=6cm]{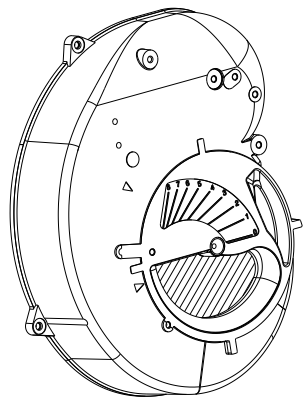}
	\caption{{\it Air intake (hatched section) for a fan housing}.}
	\label{damper}
	\end{center}
\end{figure}

\begin{figure}[ht!]
	\begin{center}
	\includegraphics[height=6cm]{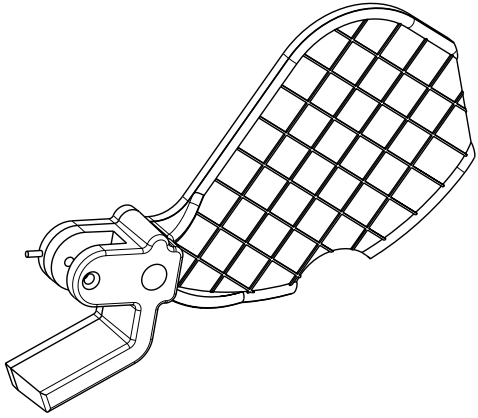}
	\caption{{\it Shutter of an air damper for low power burners}.}
	\label{damper}
	\end{center}
\end{figure}

\begin{figure}[ht!]
	\begin{center}
	\includegraphics[height=8cm]{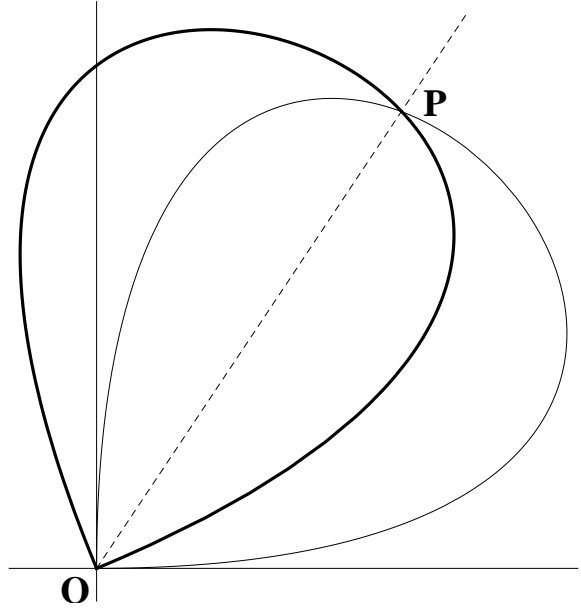}
	\caption{{\it The damper (thick line for the shutter) is partially open}.}
	\label{damperOval}
	\end{center}
\end{figure}

\noindent It could be very important to have a simple mathematical formula for a fast computing of the damper section left free by the shutter in its opening positions. Suppose that hole and shutter have the same geometrical shape and that the boundary of this shape is a smooth oval. For example, consider the oval $x^3+y^3+3(x^2y-xy+xy^2)=0$, which have $x=3(1-t)^2t, \hspace{0.2cm} y=3(1-t)t^2$, $t \in [0,1]$, as (0,0)-centered parametrization. From previous Theorem, this oval is locally algebraically squarable. The shutter can rotate with center on the origin. If the damper is partially open, the section left free is equal to the difference of the two segments cut from the ovals by the line $OP$ (Fig.\ref{damperOval}). These segments are algebraically squarable, therefore the free section, that is their difference, is algebraically squarable.\\

\noindent {\it Exercise}. Show that the straight line $OP$ is axis of simmetry for the figure composed by the two thin and thick ovals. You could draw the ovals on a paper sheet and then fold it along the line $OP$.\\

\noindent The equation of line $OP$ is $y=m_Px$, where $m_P=\frac{t_P}{1-t_P}$, with $t$ parameter of the hole boundary curve. Note that, if $\alpha$ is the angle between $x$-axis and the straight line, then for every $\alpha \in [0,\frac{\pi}{2}]$ there is a unique $t_P$ such that $x_P = x(t_P)$ and $y_P = y(t_P)$. The area of the segment of the hole (thin curve) is, using Gauss-Green,

\begin{figure}[ht!]
	\begin{center}
	\includegraphics[height=2cm]{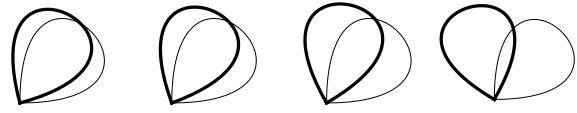}
	\caption{{\it Various positions of shutter}.}
	\label{positions}
	\end{center}
\end{figure}

\begin{eqnarray}
	\nonumber S_1 = - \left( \int_0^{t_P}y(t)x'(t)dt \hspace{0.1cm} + \hspace{0.1cm} m_P \int_{x_P}^0 xdx \right) =\\ 
	= - 3{t_P}^3 + \frac{45}{4}{t_P}^4 - \frac{63}{5}{t_P}^5 + \frac{9}{2}{t_P}^6 + \frac{1}{2} m_P x_P^2
\end{eqnarray}

\noindent The area of the total thin oval is

\begin{eqnarray}
	S = - \int_0^{1}y(t)x'(t)dt \hspace{0.1cm} = \frac{3}{20}
\end{eqnarray}

\noindent therefore the area of the free inlet section is

\begin{eqnarray}
	\nonumber S_2 = S_1 - (S - S_1) = 2S_1 - S =\\ 
	\nonumber = - 6{t_P}^3 + \frac{45}{2}{t_P}^4 - \frac{126}{5}{t_P}^5 + 9{t_P}^6 +  m_P x_P^2 - \frac{3}{20}
\end{eqnarray}

\begin{figure}[ht!]
	\begin{center}
	\includegraphics[width=8cm]{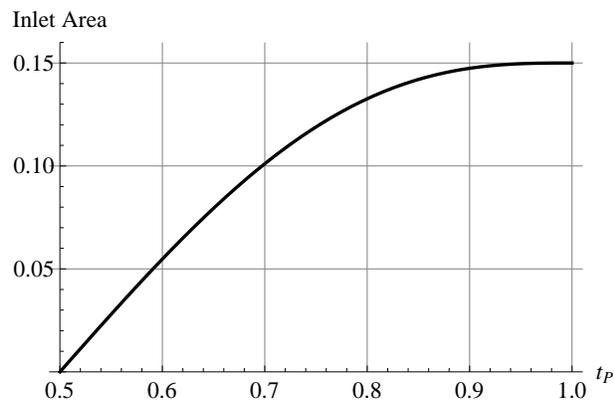}
	\caption{{\it Graph of the area $S_2$ as function of $t_P$}.}
	\label{areaGraph}
	\end{center}
\end{figure}

\noindent Note that you could express $S_2$ as a rational function only of $m_P$, $m_P \in [1,+\infty]$, variable or only of $t_P$, $t_P \in [\frac{1}{2},1]$, variable; in particular $m_P x_P^2 = 9(1-t_P)^3t_P^2$, so that $S_2$ is a polynomial function of $t_P$.

\end{document}